\documentstyle[a4,12pt]{article}

\def\mbo#1{{\mathchoice {\mbox{\normalshape\normalsize #1}}
{\mbox{\normalshape\normalsize #1}} {\mbox{\normalshape\scriptsize #1}}
{\mbox{\normalshape\tiny #1}} }}
\def\bbbr{{\mbo{I$\!$R}}} 
\def\bbbc{{\mathchoice {\setbox0=\hbox{$\displaystyle\mbo{C}$}\hbox{\hbox
to0pt{\kern0.4\wd0\vrule height0.9\ht0\hss}\box0}}
{\setbox0=\hbox{$\textstyle\mbo{C}$}\hbox{\hbox
to0pt{\kern0.4\wd0\vrule height0.9\ht0\hss}\box0}}
{\setbox0=\hbox{$\scriptstyle\mbo{C}$}\hbox{\hbox
to0pt{\kern0.4\wd0\vrule height0.9\ht0\hss}\box0}}
{\setbox0=\hbox{$\scriptscriptstyle\mbo{C}$}\hbox{\hbox
to0pt{\kern0.4\wd0\vrule height0.9\ht0\hss}\box0}}}}

\begin{document}

\titlepage

\begin{center}
{\Large{\bf On the non-commutative \\ Riemannian
geometry of $GL_{q}(n)$}}
\end{center}

\vskip1cm
\begin{center}
Y. Georgelin$^{1}$, J. Madore$^{2}$, T. Masson$^{2}$, J. Mourad$^{3}$
\end{center}
\vskip.5cm

\begin{center}
$^{1}$ {\it Division de Physique Th\'eorique\footnote{
Unit\'e de recherche des universit\'es Paris XI et Paris VI associ\'ee
au CNRS} \\
Institut de Physique Nucl\'eaire,
Univerist\'e de Paris Sud\\
F-91406 Orsay}
\end{center}
\vskip .5cm
\begin{center}
$^{2}$ {\it Laboratoire de Physique Th\'eorique et Hautes Energies
\footnote{Laboratoire associ\'e au CNRS URA D0063},\\
B\^at. 211, Universit\'e de Paris Sud,\\
F-91405 Orsay}
\end{center}
\vskip1cm

\begin{center}
$^3$ {\it Laboratoire de Mod\`eles de Physique Math\'ematique,\\
Parc de Grandmont, Universit\'e de Tours,\\
F-37200 Tours\\
and\\
D\'epartement de Physique
\footnote{address after September 1, 1995.},\\
49 av. des Genottes, Universit\'e de Cergy Pontoise,\\
BP 8428, F-95806 Cergy Pontoise\\
e-mail: celfi.phys.univ-tours.fr}
\end{center}
\vskip.5cm

\begin{center}
{\bf Abstract}
\end{center}

A recently proposed
definition of a linear connection in non-commutative geometry,
based  on a generalized permutation, is used to construct
 linear connections on $GL_{q}(n)$.
Restrictions on the generalized permutation
arising from the stability of linear connections under
involution are discussed.
Candidates for generalized permutations on $GL_q(n)$
are found. It is shown that, for a given generalized
permutation, there exists one and only one associated
linear connection. Properties of the linear
connection are discussed, in particular its
bicovariance, torsion and  commutative
limit.

\vskip.5cm
\noindent
{\bf LPTHE-Orsay 95/51}

\noindent{q-alg/9507002}

\newpage

\section{Introduction}
\setcounter{equation}{0}
Shortly after their discovery
in the context of integrable models \cite{skl,kul,dri,jim},
quantum groups were identified as interesting non-commutative
generalizations of the algebra of functions
on a  Lie-group manifold \cite{woro,man,fad}.
 Non-commutative differential
calculi \cite{con} have been proposed where the main constraint is
the bicovariance of the differential algebra
\cite{wor}.
In addition the $R$-matrix formulation \cite{fad}
played a key role in further developements \cite{ber,mal,jur,
schu,schi,sud}.

The aim of this paper is to define
  linear connections and metrics on quantum groups.
{}From the mathematical point of view,
this is a step towards the understanding of
which classical concepts can have a non-commutative
generalization.
{}From the physical point
of view, it could be a first step towards the formulation
of gravitational theories on quantum groups.
Non-commutative manifolds in fact could represent a solution
to the problem of short distance divergences of
usual quantum field theories (See {\it e.g.} \cite{maj,mad1})
and could also offer a more satisfactory description of
space-time.
In this respect, quantum groups
are an interesting toy model where qualitative differences
between the non-commutative ($q \neq 1$) and the
non-deformed ($q=1$) cases can be observed.
 In the context of the Dirac-operator-based differential
calculus of Connes, an approach to the construction of such theories
has been proposed
using the Wodzicki residue of the Dirac operator \cite{kal,kas}.
However, many interesting differential calculi,
such as those on quantum groups \cite{wor} and spaces \cite{wes},
are not defined by a Dirac operator. Here,
as was proposed in \cite{mad2,cham,keh}, we
follow the idea, which is suitable for all differential calculi,
 of a generalization  to the non-commutative
context of the usual commutative
metrics and linear connections.

A general definition of linear connections, in the context of
non-commutative geometry,
has been recently proposed for
 the derivation-based differential calculus \cite{mou,dub1} and
other differential calculi \cite{mou} in which case the construction
 relies on a generalized permutation.
 In Section 2 we fix our notation concerning quantum groups.
In Section 3 we briefly review the construction of
\cite{mou} and add some
restrictions on the generalized permutations
which arise from the requirement that
the set of covariant derivatives be
stable under complex conjugation. Section 4 is devoted to the
search for generalized permutations on $GL_{q}(n)$
which are restricted by the bicovariance condition.
We find a  two-parameter family of generalized permutations.
In Section 5
we prove that for a given generalized permutation there exists only one
linear connection. Properties of this linear connection
are studied, in particular its bicovariance, torsion and curvature.
Finally,  we examine the commutative limit of
our linear connections.
We show that  the limit of one class of
these when $q \rightarrow 1$ corresponds to
left and right invariant linear connections
  on $GL(n)$.  We collect our conclusions in Section 6.

\section{Quantum groups and their bicovariant differential calculi}
\setcounter{equation}{0}
The  quantum group  ${\rm Fun}(GL_{q}(n))$
 is a  Hopf algebra (${\cal A}, \Delta,
\epsilon, \kappa$) generated, as an
algebra, by the identity and $T_{i}{}^{j}$, $i,j=1,\dots n$.
An exchange of the order of the generators, while maintaining the classical
Poincar\'e series,
 is
obtained
by the RTT relation \cite {fad}:
\begin{equation}
R T_{1} T_{2} =T_{1} T_{2}  R.\label{commu}
\end{equation}
Here,  $R$ is the $R$-matrix, which is an element
of $M_{n}(\bbbc)\otimes M_{n}(\bbbc)$
obeying the
Yang-Baxter relation
\begin{equation}
(1\otimes R) \ (R\otimes 1) \ (1 \otimes R)\ =\ (R\otimes 1) \ (1\otimes R)
\ (R\otimes 1).
\end{equation}
The $R$ matrix of $GL_{q}(n)$ is given by \cite{jim}
\begin{equation}
\label{matriceR}
R=q\sum_{i}E_{ii}\otimes E_{ii}+\sum_{i\neq j}E_{ij}\otimes E_{ji}
+\lambda\sum_{i<j}E_{ii}\otimes E_{jj},
\end{equation}
where $\lambda=q-q^{-1}$.
It satisfies the Hecke condition
\begin{equation}
\label{hec}
 (R-q)(R+{1\over{q}})=0.
\end{equation}

  The differential calculus on the quantum
group is considerably restricted by the bicovariance condition \cite{wor}.
This means that there exist a right and left
 coaction of $\cal A$
on $\Omega^{1}$, the space of 1-forms, such that
\begin{eqnarray}
\Delta_{L}(a\ db)=\Delta(a)\ (1\otimes d)\ \Delta(b), \\
\Delta_{R}(a \ db)=\Delta(a)\ (d\otimes 1)\ \Delta(b), \\
(1\otimes \Delta_{R})\ \Delta_{L} \ =\ (\Delta_{L}\otimes 1)\ \Delta_{R}.
\end{eqnarray}

Under some restrictions on $q$ and the
assumption that $\Omega^{1}$ be generated as
a left-module by $dT^{i}_{j}$, bicovariant differential calculi
have been  classified \cite{sch} and shown to be obtained by
the constructive method of Jurco \cite{jur}.

For such differential calculi $\Omega^{1}$ is generated
as a left (or right) module by left invariant 1-forms
$\omega^{i}_{j}$ ($\Delta_{L}(\omega^{i}_{j})=1\otimes\omega^{i}_{j}$):
\begin{equation}
\omega^{i}_{j}\ =\ \kappa(T^{i}_{k})\ dT^{k}_{j}.
\end{equation}
The differential algebra is entirely
 characterized by the commutation relations
\begin{equation}
\omega^{i}_{j}a=(1\otimes f^{i}_{j}{}^{k}_{l})\ \Delta(a)\omega^{l}_{k},
\end{equation}
where $f^{i}_{j}{}^{k}_{l}$ are linear functionnals representing
the algebra $\cal A$:
$$f^{i}_{j}{}^{k}_{l}(1)=\delta^{i}_{l}\delta_{j}^{k},\quad
f^{i}_{j}{}^{k}_{l}(ab)=
f^{i}_{j}{}^{m}_{n}(a)f^{n}_{m}{}^{k}_{l}(b).$$

They can be explicitly determined in terms of the $R$-matrix
and some parameters  \cite{sch}. Here we shall make the choice
\begin{equation}
\label{cal}
f^{jk}_{il}(T^{m}_{n})=
 (R^{-1})^{sk}_{in}\ (R^{-1})^{jm}_{sl}
\end{equation}
for the differential calculus. In the limit $q\rightarrow 1$,
this differential calculus reduces to the usual one on $GL(n)$.
It has been considered in Ref. \cite{mal,sud,schi,zum,schu}.
In this case,
the commutation relations are often written in the form
\begin{equation}
\label{commutdt}
 T_{1}dT_{2}=RdT_{1}T_{2}R.
\end{equation}

The space of 2-forms is constructed as the image
of $\Omega^{1}\otimes_{\cal A}\Omega^{1}$ under the ``multiplication''
map $\pi$:
\begin{eqnarray}
\pi :\Omega^{1}\otimes_{\cal A}\Omega^{1} &\rightarrow& \Omega^{1}\otimes_{\cal
A}\Omega^{1}\\
\pi &=& 1-\Lambda,
\end{eqnarray}
where $\Lambda$ is a bimodule automorphism, obeying
the Yang-Baxter equation, which generalizes
the permutation map of the commutative case.
Let $\eta^{i}_{j}$ be right invariant 1-forms:
\begin{equation}
\eta^{i}_{j} =T^i_l \omega^l_k \kappa(T^k_j).
\end{equation}
Then $\Lambda$ is determined  by \cite{wor}
\begin{equation}
\Lambda(\omega^{i}_{j}\otimes\eta^{k}_{l})\ =\
\eta^{k}_{l}\otimes\omega^{i}_{j}.
\end{equation}
When applied to $\omega^{i}_{j}\otimes\omega^{k}_{l}$,
one can  show that
\begin{eqnarray}
\Lambda(\omega^{i}_{j}\otimes\omega^{k}_{l})=\Lambda^{i}_{j}{}^{k}_{l}{}
\ _{m}^{n}{}_{q}^{p} \ \omega^{m}_{n}\otimes\omega^{q}_{p}\\
\Lambda^{i}_{j}{}^{k}_{l}{}
\ _{m}^{n}{}_{q}^{p} \ =\ f^{ip}_{jq}( \kappa( T^k_m) T^n_l ).
\end{eqnarray}
When applied to $dT_{1}\otimes dT_{2}$, the  map
$\Lambda$ yields
\begin{equation}
\label{lam}
\Lambda(dT_{1}\otimes dT_{2})=R dT_{1}\otimes dT_{2} R^{-1}.
\end{equation}
The Hecke relation for the $R$ matrix (\ref{hec}), combined
with the previous equation, yields
the following characteristic equation for $\Lambda$:
\begin{equation}
(\Lambda-1)(\Lambda+q^{2})(\Lambda+q^{-2})=0.
\end{equation}

Higher order forms can be constructed in a
similar way using the map $\Lambda$ \cite{wor}.
The exterior derivative is defined with the help
of the right and left-invariant 1-form $\theta$:
\begin{equation}
\label{theta}
\theta=-{q^{2n+1}\over {\lambda}} \sum_{i} q^{-2i}\omega^{i}_{i},
\end{equation}
by
\begin{equation}
d\omega=[\theta,\omega],
\end{equation}
where $[,]$ is the graded commutator and the product is in $\Omega$.

For real values of $q$ or for $|q|=1$, an involution
may be defined on $GL_{q}(n)$ reducing it respectively
to $U_{q}(n)$ or to $GL_{q}(n,\bbbr)$. Setting the $q$-determinant
equal to one gives rise to $SU_{q}(n)$ and $SL_{q}(n,\bbbr)$ \cite{fad}.
The bicovariant differential calculus on these reductions
is characterized either by a larger set of 1-forms than the
classical case \cite{car} or by a modified Leibniz rule \cite{fadd}.

\section{Linear connections in non-commutative geometry}
\setcounter{equation}{0}
In this section we collect the main definitions and
results concerning the general construction of linear
connections as proposed in \cite{mou}. We add some new restrictions
on the generalized permutation by imposing the stability of
the set of covariant derivatives under complex conjugation.
In the following
$\cal A$ is a unital associative algebra over $\bbbc$
equipped with the differential calculus $(\Omega,d)$.

\bigskip
\noindent
{\bf Definition 3.1: }
{\it Let $\pi$ be the multiplication in $\Omega$.
A generalized permutation, $\sigma$, is a bimodule automorphism
of $\Omega^1 \otimes_{\cal A}\Omega^{1}$
satisfying
\begin{equation}
\label{defsigma}
\pi\circ\sigma=-\pi.
\end{equation}
A generalized flip, $\tau$, is defined as a generalized
permutation satisfying $\tau^2=1$.
}

\bigskip
\noindent
{\bf Remarks:}

1-  Note that $\sigma=-1$ is a generalized flip.

2- When the algebra ${\cal A}$
is the algebra of $C^\infty$ functions on a manifold $M$,
the permutation
\begin{equation}
\tau(\omega\otimes\omega')=\omega'\otimes\omega
\end{equation}
is a generalized flip.

3- When $\Omega^2$ is realized as a subspace of $\Omega^{1}
\otimes\Omega^1$ with an imbedding $i$ verifying
$\pi\circ i=1_{\Omega^2}$ then
\begin{equation}
\label{can}
1\ -\ 2\ i \circ \ \pi
\end{equation}
is a generalized flip. The generalized flip of the derivation-based
differential calculus proposed in Ref.\cite{dub1,mou} is of this
form, as are the generalized flips of Ref.\cite{mmm,keh}.

4- If $\sigma$ is a generalized permutation then so is $\sigma^{-1}$
as well as $\sigma^{2n+1}$ for an arbitrary integer $n$.

5- If $\sigma$ and $\sigma'$ are two generalized permutations
then so is $\mu(\sigma + 1) + \mu'(\sigma' + 1) - 1$.
The $\sigma+1$ form  a linear space.

\bigskip
\noindent
{\bf Definition 3.2: }
{\it A linear connection associated
to a generalized permutation, $\sigma$,
is a linear map, $\nabla^{\sigma}$,
 from $\Omega^1$ to $\Omega^{1}\otimes_{\cal A}
\Omega^{1}$ satisfying the two Leibniz rules
\begin{eqnarray}
\label{defsigmaleft}
\nabla^{\sigma}(a\omega) &=& da\otimes\omega+a \ \nabla^{\sigma}\omega,\\
\label{defsigmaright}
\nabla^{\sigma}(\omega a) &=& \sigma(\omega\otimes da)+\nabla^{\sigma}\omega
\ a,
\end{eqnarray}
for any $a \in  {\cal A}$ and any $\omega  \in  \Omega^1$.
}

\bigskip
\noindent
{\bf Remarks:}

1- When the algebra $\cal A$ is the commutative
algebra of smooth functions on a manifold
 the only possible
linear connections are those associated to the permutation (3.2).

2- If $\sigma$ and $\sigma'$ are
two generalized permutations then
$\nabla^{\sigma}-\nabla^{\sigma'}$
is a left-module homomorphism.

3- If $\nabla$ and $\nabla'$ are two linear connections associated to
the same generalized permutation, then their difference is
a bimodule homomorphism.

\bigskip
The preceding definition of the linear connection
has the advantage of allowing
an extension to the tensor product over $\cal A$
of several copies of $\Omega^1$. This is formulated in:

\noindent
{\bf Proposition 3.3: } {\it A linear connection
associated to a generalized permutation $\sigma$
admits a unique extension as a linear map from
$\underbrace{\Omega^1 \otimes_{\cal A}\dots\otimes_{\cal A}\Omega^1}
_{\rm {s \ times}}$ to $\underbrace{
\Omega^1\otimes_{\cal A}\dots \otimes_{\cal A}\Omega^1}_{\rm {s+1 \ times}}$
of the form
\begin{equation}
\nabla^{\sigma}(\omega\otimes\omega')=\nabla^{\sigma}(\omega)
\otimes\omega'+\sigma_s(\omega\otimes\nabla^{\sigma}\omega'),
\end{equation}
for any $\omega \in \Omega^{1}$ and any $\omega' \in
\underbrace{\Omega^{1}\otimes_{\cal A}\dots\otimes_{\cal A}\Omega^{1}}
_{\rm {s-1 \ times}}$,
with $\sigma_{s}$ an automorphism
of $\underbrace{\Omega^{1}\otimes_{\cal A}\dots\otimes_{\cal A}
\Omega^{1}}_{\rm {s+1 \ times}}$. The unique $\sigma_{s}$
 is given by
\begin{equation}
\sigma_{s}=\sigma\otimes
\underbrace{1\otimes\dots\otimes 1}_{\rm {s-1 \  times}}.
\end{equation}
}

\noindent
{\it Proof.}
The proof can be carried out by induction. For $s=2$ an
identification of $\nabla^{\sigma}(\omega f\otimes\omega')$
with $\nabla^{\sigma}(\omega\otimes f\omega')$, where
 $f$ is an arbitrary element of $\cal A$ and $\omega$ and $\omega'$
are 1-forms, gives $\sigma_{2}=\sigma\otimes 1$;
so the proposition is true for $s=2$.
Suppose it is true to order $s-1$ and let $\omega'$ be
an element of the tensor product of $s-1$ copies
of $\Omega^{1}$ then, by the induction hypothesis,
\begin{equation}
\nabla^{\sigma}f\omega'=df\otimes\omega'+f\nabla^{\sigma}\omega'.
\end{equation}
The identification of $\nabla^{\sigma}(\omega\otimes f\omega')$
with $\nabla^{\sigma}(\omega f\otimes
\omega')$ where $\omega$ is an element of
$\Omega^{1}$ completes the Proof. $\clubsuit$

\bigskip
Suppose that $\cal A$ is an algebra over $\bbbc$ equipped with
an involution *. Then $\Omega^{1}$ carries a natural
involution defined by $(bda)^{*}=(da^{*})b^{*}$.
The involution on $\Omega^{1}\otimes_{\cal A}\Omega^{1}$
is not {\it a priori} determined. In fact, if we define
the anti-homomorphism $\alpha$ by
\begin{equation}
\label{alpha}
\alpha(\omega\otimes\omega')=\omega'^{*}\otimes\omega^{*},
\end{equation}
and if $\phi$ is an automorphism of
$\Omega^{1}\otimes_{\cal A}\Omega^{1}$ such that $(\phi\circ\alpha)^{2}=1$,
then $\phi\circ\alpha$ defines an involution on
$\Omega^{1}\otimes_{\cal A}\Omega^{1}$.
We would like to define an involution on
$\Omega^{1}\otimes_{\cal A}\Omega^{1}$ which in the commutative limit
reduces to $(\omega\otimes\omega')^{*}=
\tau(\omega'^{*}\otimes\omega^{*})$,
where $\tau$ is the usual permutation operator, and which
allows us to define the complex
conjugate of a  linear connection, as in the commutative case,
by
\begin{equation}
\label{defconjlinconn}
\overline{\nabla}^{\sigma} \omega=\Big(\nabla^{\sigma}(\omega^{*})\Big)^{*}.
\end{equation}
The requirement that $\overline{\nabla}^{\sigma}$ be a linear connection
imposes constraints on the involution on
$\Omega^{1}\otimes_{\cal A}\Omega^{1}$ and on
the generalized permutation, $\sigma$:

\bigskip
\noindent
{\bf Proposition 3.4: }
{\it Suppose that $\cal A$ is equipped with
an involution $^{*}$, then the following  assertions are equivalent

1- The map ${\overline{\nabla}}^{\sigma}$
defined in (\ref{defconjlinconn}) is a linear connection.

2- the generalized permutation,
$\sigma$, verifies
\begin{equation}
\label{conditioninv}
(\sigma\circ\alpha)^{2}=1
\end{equation}
and
the involution on $\Omega^{1}\otimes_{\cal A}\Omega^{1}$
is given by
\begin{equation}
\label{sigmainvolution}
(\omega\otimes\omega')^{*}=\sigma(\omega'^{*}\otimes\omega^{*}).
\end{equation}
}

\noindent
{\it Proof. } $2 \Rightarrow 1$ is a direct calculation.
We prove $1\Rightarrow 2$.
Calculate, with the aid of Equation (\ref{defconjlinconn}),
 ${\overline{\nabla}}^{\sigma}(\omega a)$:
\begin{eqnarray}
{\overline{\nabla}}^{\sigma}(\omega a) &=&
\Big(\nabla^{\sigma}(a^{*}\omega^{*})\Big)^{*}
\nonumber \\
 &=& \Big(da^{*}\otimes\omega^{*}\Big)^{*}+
\Big(\nabla^{\sigma}\omega^{*}\Big)^{*} a \nonumber \\
 &=&\Big(da^{*}\otimes\omega^{*}\Big)^{*}+ ({\overline{\nabla}^{\sigma}
\omega}) a.
\label{un}
\end{eqnarray}
If the map $\overline{\nabla}^{\sigma}$ is a covariant derivative
then there exists a generalized permutation, $\phi$, such that
\begin{equation}
\label{de}
 {\overline{\nabla}}^{\sigma}(\omega a)=
 \phi(\omega\otimes da)+\overline{\nabla}\omega \ a.
\end{equation}
Comparing the two equations (\ref{un}) and (\ref{de}) we obtain
\begin{equation}
\Big(da^{*}\otimes\omega^{*}\Big)^{*}=\phi(\omega\otimes da).
\end{equation}
This equation is valid for arbitrary $a$ and $\omega$
so the involution in $\Omega^{1}\otimes_{\cal A}\Omega^{1}$ verifies:
\begin{equation}
\label{invo}
\Big(\omega'\otimes\omega\Big)^{*}=\phi(\omega^{*}\otimes \omega'^{*}).
\end{equation}
The involution property, $**=1$, gives $(\phi\circ\alpha)^{2}=1$.
It remains to prove that $\phi=\sigma$.
In order to do this, calculate, using Equation (\ref{defconjlinconn}),
$\overline{\nabla}^{\sigma}(a\omega)$:
\begin{equation}
\label{un1}
\overline{\nabla}^{\sigma}(a\omega) =
a\Big(\nabla^{\sigma}(\omega^{*})\Big)^{*}
+\Big(\sigma(\omega^{*}\otimes da^{*})\Big)^{*}.
\end{equation}
Since $\overline{\nabla}^{\sigma}$ is a linear connection
we have
\begin{equation}
\label{de2}
\overline{\nabla}^{\sigma}(a\omega) =
a {\overline{\nabla}}^{\sigma}\omega +da\otimes\omega.
\end{equation}
Comparing equation (\ref{un1}) and (\ref{de2}) we get
\begin{equation}
da\otimes\omega=\Big(\sigma(\omega^{*}\otimes da^{*})\Big)^{*}.
\end{equation}
This equation is valid for arbitrary $a$ and $\omega$, so we have
\begin{equation}
\label{invol}
\omega'\otimes\omega=\Big(\sigma(\omega^{*}\otimes \omega'^{*})\Big)^{*}.
\end{equation}
Comparing  equations (\ref{invo}) and (\ref{invol})
leads to the equality of $\phi$ and $\sigma$. $\clubsuit$

\bigskip
\noindent
{\bf Definition 3.5: } {\it For a given involution $*$
on $\Omega^{1}\otimes_{\cal A}\Omega^{1}$,
a generalized permutation $\sigma$ is defined to be real
if it satifies the following property:
\begin{equation}
\label{defreelsigma}
\sigma \circ * = * \circ \sigma
\end{equation}
on $\Omega^{1}\otimes_{\cal A}\Omega^{1}$.
}

\bigskip
Now, if one wants to find an involution $*$
on $\Omega^{1}\otimes_{\cal A}\Omega^{1}$
such that ${\overline{\nabla}}^{\sigma}$ is a linear connection,
then one should
 take, according
to Proposition 3.4, (\ref{sigmainvolution}) as a definition of $*$.
The condition for this to be possible is $(\sigma\circ\alpha)^{2}=1$.
If one further demands that $\sigma$  be real
then one has to use the following:

\noindent
{\bf Proposition 3.6: } {\it
Suppose that the generalized permutation,
$\sigma$, verifies Equation (\ref{conditioninv}) and that
the involution on
$\Omega^{1}\otimes_{\cal A}\Omega^{1}$ is given by $\sigma\circ\alpha$
then $\sigma$ is real iff it is a generalized flip.
}

\noindent
{\it Proof. }
 The reality condition reads
\begin{equation}
\sigma\circ\sigma\circ\alpha=\sigma\circ\alpha\circ\sigma.
\end{equation}
Since $\sigma$ is an automorphism, this equation leads to
\begin{equation}
\sigma\circ\alpha=\alpha\circ\sigma.
\end{equation}
This relation in $(\sigma\circ\alpha)^2 = 1$ gives $\sigma^{2}=1$. $\clubsuit$

\bigskip
The definition of the complex conjugate of a linear connection
allows the following:

\noindent
{\bf Definition 3.7: }
{\it A real linear connection
associated to a generalized permutation $\sigma$ is defined by
${\overline{\nabla}}^{\sigma}=\nabla^\sigma$.}

\bigskip
\noindent
{\bf Remark:}

The involution on $\Omega^{1}\otimes_{\cal A}\Omega^{1}$
defined above induces
an involution on $\Omega^{2}$ by $(\omega\wedge\omega')^{*}=
\pi((\omega\otimes\omega')^{*})=-\omega'^{*}\wedge\omega^{*}$.
This is due to the property (\ref{defsigma}).

\bigskip
\noindent
{\bf Definition 3.8: } {\it The torsion $T$ of a linear connection
$\nabla^\sigma$, is defined as the linear map
from $\Omega^{1}$ to $\Omega^{2}$ given by
\begin{equation}
\label{tor}
T=d-\pi\circ\nabla^\sigma.
\end{equation}
}

\bigskip
\noindent
{\bf Proposition 3.9:}
{\it The torsion map is a bimodule homomorphism.}

\noindent
{\it Proof.} It is an immediate consequence of
the condition (\ref{defsigma}). $\clubsuit$

\bigskip
\noindent
{\bf Definition 3.10: } {\it The curvature ${\cal R}$
of a linear connection $\nabla^\sigma$
is defined as the linear map from
$\Omega^{1}$ to $\Omega^{2}\otimes_{\cal A}\Omega^1$ given by
\begin{equation}
{\cal R} = \left(
(T\otimes 1) + (\pi\otimes 1)\nabla^\sigma
 \right) \nabla^\sigma
\end{equation}
}

\noindent
{\bf Proposition 3.11:} {\it The curvature is
a left-module homomorphism.}

\noindent
{\it Proof.} A straightforward calculation. $\clubsuit$

\bigskip
\noindent
{\bf Definition 3.12: } {\it A metric $g$ is defined as
an element of $\Omega^{1}\otimes_{\cal A}\Omega^{1}$
satisfying
\begin{equation}
\pi(g)=0.
\end{equation}
If $\Omega^{1}\otimes_{\cal A}\Omega^{1}$ is equipped
with an involution, a real metric is defined by $g^{*}=g$.}

The definition of a non-degenerate metric
requires some
more structure on the algebra $\cal A$.
This structure must guarantee
that the dimension of $\Omega^{1}$ as a left module be well defined.
For example if $\cal A$ is a Hopf
algebra then it is well known that
this is so (See {\it e.g.} \cite{wor}).
If it exists, let $\omega^{a},\ a=1,\dots N$, be a free basis of $\Omega^{1}$
as a left module then
a metric can be written uniquely in the form
\begin{equation}
g=g_{ab}\omega^{a}\otimes\omega^{b},
\end{equation}
with $g_{ab}\in {\cal A}$. We will call a metric
nondegenarate if the  matrix whose elements are $g_{ab}$
is invertible.

\bigskip
\noindent
{\bf Definition 3.13: }
{\it A metric $g$ and  a linear connection
$\nabla^\sigma$ are said to be compatible
if the condition $\nabla^{\sigma} g=0$ is satisfied.}

\section{Determination of $\sigma$ on $GL_{q}(n)$}
\setcounter{equation}{0}

In addition to the previous requirements on $\sigma$,
it is natural, in the context of quantum groups, to
add the requirement of bicovariance:

\noindent
{\bf Definition 4.1:} {\it
A generalized permutation, $\sigma$, is called bicovariant
iff:
\begin{eqnarray}
\label{bicov}
(1\otimes\sigma)\Delta_{L} &=& \Delta_{L}\ \sigma, \nonumber \\
(\sigma\otimes 1) \Delta_{R} &=& \Delta_{R} \ \sigma.
\end{eqnarray}
}

Following the $R$-matrix technique, that is the determination
of all unknown maps from the $R$-matrix and $q$,
we will determine the candidates
for the map $\sigma$ in terms of $R$.
We recall from (2.13) and (3.1)
that the generalized permutation $\sigma$ is an automorphism
of $\Omega^{1}\otimes_{\cal A}\Omega^{1}$ verifying
\begin{equation}
(1-\Lambda)\circ(\sigma+1)=0,
\end{equation}
the bicovariance requirements
(\ref{bicov}), and when ${\cal A}$ is equipped with
an involution, that is for real $q$ and for $|q|=1$, the involution
property ({\ref{conditioninv}).

\bigskip

In order to find  candidates for $\sigma$, we
shall prove the following Proposition, which, in its first part,
is a generalization of Proposition 3.1 of \cite{wor}:

\noindent
{\bf Proposition 4.2: }
{\it
Let $\alpha_{ij}, \ i,j=0,1,$ be complex numbers.

1-There exists
a unique bimodule homomorphism, $\Phi$, of
$\Omega^{1}\otimes_{\cal{A}}\Omega^{1}$
such that
\begin{equation}
\label{defbi}
\Phi(dT_{1}\otimes dT_{2})=\sum_{i,j}\alpha_{ij}
R^{i}\ dT_{1}\otimes dT_{2} \ R^{j}.
\end{equation}

Moreover,

2- The map $\Phi$ is bicovariant.

3- The map $\Phi$
 is a generalized permutation iff
\begin{eqnarray}
&&\alpha_{01}-\alpha_{10}=0, \nonumber \\
&&\alpha_{00}+\lambda\alpha_{10}-\alpha_{11}=-1, \label{gp}
\end{eqnarray}
where, we recall, $\lambda = q-q^{-1}$.

In this case, $\Phi$ obeys the characteristic equation
\begin{eqnarray}
&(\Phi+1)(\Phi-\lambda_{1})(\Phi-\lambda_{2})=0,\\
&\lambda_{1}=-1+\alpha_{10}(q+q^{-1})+\alpha_{11}(1+q^{2}),\\
&\lambda_{2}=-1-\alpha_{10}(q+q^{-1})+\alpha_{11}(1+q^{-2}).
\end{eqnarray}
}

\noindent
{\it Proof.}
An element $\nu$ of $\Omega^{1}\otimes_{\cal{A}}\Omega^1$
can be written in a unique
way as
\begin{equation}
\label{deco}
\nu=\sum a^{ij}_{kl}dT^{k}_{i}\otimes dT^{l}_{j}=
Tr\left(a\left( dT_{1}\otimes dT_{2}\right)\right),
\end{equation}
where $a \in M_{n}({\cal A})\otimes M_{n}({\cal A})$.
This is a consequence of the fact that the $dT$ generate
$\Omega^1$ as a left-module.
The action of $\Phi$ on $\nu$ is defined by
\begin{equation}
\Phi(\nu)=Tr(a\alpha_{ij}R^{i}dT_{1}\otimes dT_{2}R^{j}).
\end{equation}
It clearly satisfies (\ref{defbi}).
It remains to check that $\Phi$ defined in this way
is  a bimodule homomorphism. The left-module
homomorphism property is assured by construction.
To check the right-module homomorphism property it
suffices to verify that
\begin{equation}
\label{bi}
\Phi(dT_{1}\otimes dT_{2}T_{3})=
\Phi(dT_{1}\otimes dT_{2})T_{3}.
\end{equation}
This is so because the $T$ generate the algebra.
The left-hand side of Equation (\ref{bi}) can
be written, after successive use of Equation (\ref{commutdt}), as
\begin{equation}
\label{bii}
\Phi(dT_1\otimes R^{-1}_{23}T_{2} dT_{3} R^{-1}_{23})=
R^{-1}_{23}R^{-1}_{12}\Phi(T_{1}dT_{2}\otimes dT_{3})
R^{-1}_{12}R^{-1}_{23},
\end{equation}
here the subscripts of the $R$-matrix denote
the two spaces on which it acts.  Next,
we use the left-module property to write
the right-hand side of Equation (\ref{bii})
as
\begin{equation}
\label{gauche}
R^{-1}_{23}R^{-1}_{12}T_{1}\Phi(dT_{2}\otimes dT_{3})
R^{-1}_{12}R^{-1}_{23}=\alpha_{ij}R^{-1}_{23}R^{-1}_{12}
R^{i}_{23}T_{1}dT_{2}\otimes dT_{3}R^{j}_{23}
R^{-1}_{12}R^{-1}_{23}.
\end{equation}
The right-hand side of Equation (\ref{bi})
can be written as
\begin{equation}
\alpha_{ij}R^{i}_{12}dT_{1}\otimes dT_{2} R^{j}_{12} T_{3}=\alpha_{ij}
R^{i}_{12}dT_{1}\otimes dT_{2}T_{3}R^{j}_{12}.
\end{equation}
The commutation relations (\ref{commutdt}) allow us to
write this term
as
\begin{equation}
\label{biii}
\alpha_{ij}R^{i}_{12}dT_{1}\otimes R^{-1}_{23}T_{2
}  dT_{3}R^{-1}_{23}R^{j}_{12}
=\alpha_{ij} R^{i}_{12}R^{-1}_{23}R^{-1}_{12}
T_{1} dT_{2} \otimes dT_{3} R^{-1}_{12}
 R^{-1}_{23}R^{j}_{12}.
\end{equation}
As a consequence of the Yang-Baxter equation we have
\begin{eqnarray}
 R^{i}_{12}R^{-1}_{23}R^{-1}_{12}&=R^{-1}_{23}R^{-1}_{12}R^{i}_{23},\cr
 R^{i}_{23}R^{-1}_{12}R^{-1}_{23}&=R^{-1}_{12}R^{-1}_{23}R^{i}_{12}.
\end{eqnarray}
The right hand sides of equations (\ref{gauche}) and (\ref{biii})
are thus equal. This proves the first point of the Proposition.

In order to prove the bicovariance
of $\Phi$, it suffices to prove that
\begin{eqnarray}
\label{bic1}
\Delta_{L}\Phi(dT_{1}\otimes dT_{2})=(1 \otimes \Phi)
\Delta_{L}(dT_{1}\otimes dT_{2}),\\
\label{bic2}
\Delta_{R}\Phi(dT_{1}\otimes dT_{2})= (\Phi \otimes 1)
\Delta_{R}(dT_{1}\otimes dT_{2}).
\end{eqnarray}
This is due to the fact
that $dT_{1}\otimes dT_{2}$ generate $\Omega^{1}\otimes_{\cal A}\Omega^{1}$
as a left module. Using Equation (\ref{defbi}) and
\begin{eqnarray}
\label{coal}
\Delta_{L}(dT_{1}\otimes dT_{2})=
T_{1}T_{2}\otimes dT_{1}\otimes dT_{2},\\
\label{coag}
\Delta_{R}(dT_{1}\otimes dT_{2})=dT_{1}\otimes dT_{2}\otimes T_{1}T_{2},
\end{eqnarray}
equations (\ref{bic1}) and (\ref{bic2}) can be written as
\begin{eqnarray}
\alpha_{ij}R^{i}T_{1}T_{2}\otimes dT_{1}\otimes dT_{2}R^{j}
=\alpha_{ij}T_{1}T_{2}\otimes R^{i} dT_{1}\otimes dT_{2} R^{j},\\
\alpha_{ij}R^{i} dT_{1}\otimes dT_{2}\otimes T_{1}T_{2}R^{j}
=\alpha_{ij}R^{i}dT_{1}\otimes dT_{2} R^{j} \otimes T_{1}T_{2}.
\end{eqnarray}
These equations are true due to the
commutation relations (\ref{commu}). This proves point 2 of
the Proposition.

Point 3 is a straightforward calculation
using equations (\ref{lam}) and the Hecke condition (\ref{hec}).
$\clubsuit$
\bigskip

Proposition 4.2 gives us
 a two-parameter family of bicovariant
generalized permutations. We turn to examine some of their
properties.
First, note that the maps $\Phi$ have the same eigenspaces
even though their eigenvalues might be different. In fact,
if we introduce the projectors
\begin{eqnarray}
{\Pi}_{1}(dT_{1}\otimes dT_{2})=&&
\hat P_{q} dT_{1}\otimes dT_{2} \hat P_{q},\nonumber\\
{\Pi}_{2} (dT_{1}\otimes dT_{2})=&&
\hat P_{-q^{-1}} dT_{1}\otimes dT_{2} \hat P_{-q^{-1}},\nonumber
\\
{\Pi}_{3}(dT_{1}\otimes dT_{2})=&&
\hat P_{-q^{-1}} dT_{1}\otimes dT_{2}\hat P_{q},
\nonumber\\
{\Pi}_{4}(dT_{1}\otimes dT_{2})=&&
\hat P_{q}dT_{1}\otimes dT_{2}\hat P_{-q^{-1}},
\end{eqnarray}
with
\begin{equation}
\hat P_{q}={{R+q^{-1}} \over {q+q^{-1}}},\quad
\hat P_{-q^{-1}}={{q-R} \over {q+q^{-1}}},
\end{equation}
then the generalized permutation $\Phi$
can be written as
\begin{equation}
\label{fi}
\sigma_{\lambda_{1},\lambda_{2}}
=\lambda_{1}\Pi_{1}+\lambda_{2}\Pi_{2}-\Pi_{3}-\Pi_{4},
\end{equation}
and the expression for  $\Lambda$
is
\begin{equation}
\Lambda=\Pi_{1}+\Pi_{2}-q^{2}\Pi_{3}-q^{-2}\Pi_{4}.
\end{equation}

In the commutative limit $\Pi_{1}+\Pi_{2}$ tends to the projector
onto symmetric elements of $\Omega^{1}\otimes_{\cal A}\Omega^{1}$ and
$\Pi_{3}+\Pi_{4}$ to the projector onto antisymmetric elements.
The multiplication map $\pi$ may be expressed
in terms of these projections as
\begin{equation}
\pi=(1+q^{2})\Pi_{3}+(1+q^{-2})\Pi_{4}.
\end{equation}
So $\Omega^{2}$ can be identified with the projection
of $\Omega^{1}\otimes_{\cal A}\Omega^{1}$
\begin{equation}
\Omega^{2}=(\Pi_{3}+\Pi_{4})\Omega^{1}\otimes_{\cal A}\Omega^{1}.
\end{equation}
An imbedding $i$ of $\Omega^{2}$ in
$\Omega^{1}\otimes_{\cal A}\Omega^{1}$
verifying $\pi\circ i=1_{\Omega^{2}}$ exists and is given by
\begin{equation}
i={1 \over {1+q^{2}}} \Pi_{3}+{1 \over {1+q^{-2}}} \Pi_{4}.
\end{equation}
With the aid of this imbedding we obtain the
expression (\ref{can})
for $\sigma$:
\begin{equation}
\sigma_{\Lambda}\ =\ 1\ -\ 2i\circ\pi\ =\ -1+\ 2\ (\Pi_{1}+\Pi_{2}).
\end{equation}
Note that this $\sigma$ verifies $\sigma^{2}=1$;
it is equal to $-1$ on $\Omega^{2}$ and to $+1$ on $(\Pi_{1}+\Pi_{2})
(\Omega^{1}\otimes_{\cal A}\Omega^{1})$. It corresponds to
$\lambda_{1}=\lambda_{2}=1$ in equation (\ref{fi}).

 Another simple  solution to equations (\ref{gp}), which in addition
 obeys the Yang-Baxter equation, is given by $\sigma_{q^{-2},q^{2}}$,
\begin{equation}
\label{sigmaR}
\sigma_R(dT_1\otimes dT_2) = R^{-1}dT_1\otimes dT_2 R^{-1}.
\end{equation}
This $\sigma$ is to be compared
with the $\sigma$ found in Ref. \cite{dub} for the quantum plane.
Indeed, it could be obtained in the same
way from the differential calculus (See Lemma 5.13).

\bigskip
We turn now to
  consider involutions  for
$|q| = 1$.
Then one can consider the involution $(T^i_j)^{*} = T^i_j$ on $GL_q(n)$.
This involution is compatible with
the relations on the algebra because, for $|q| = 1$
and $R$ given by (\ref{matriceR}), one has
\begin{equation}
\label{rbar}
\overline{R}^{ij}_{kl} = (R^{-1})^{ji}_{lk}
\end{equation}
where $\overline{R}$ is the complex conjugate of $R$.
In this case, the quantum group is $GL_q(n, \bbbr)$.

\bigskip
\noindent
{\bf Proposition 4.3:}
{\it Let $|q|=1$.
A generalized permutation, $\sigma_{\lambda_{1},\lambda_{2}}$,
defines an involution iff
\begin{equation}
|\lambda_{1}|\ =\ |\lambda_{2}|\ =\ 1.
\end{equation}
}

\noindent
{\it Proof.} For $|q|=1$ relations $\ref{rbar}$
imply that
\begin{equation}
\alpha\circ\Pi_{i}\circ\alpha=\Pi_{i},\ i=1,2,3,4.
\end{equation}
So, the condition $\ref{conditioninv}$ reads
\begin{equation}
|\lambda_{1}|^{2}\Pi_{1}+|\lambda_{2}|^{2}\Pi_{2}
+\Pi_{3}+\Pi_{4}=1,
\end{equation}
which completes the Proof. $\clubsuit$

\noindent
{\bf Remark:}

The previously defined $\sigma_{\Lambda}$ and $\sigma_{R}$
satisfy equations (\ref{conditioninv}).

\section{Linear connections on $GL_{q}(n)$}
\setcounter{equation}{0}

In this section we determine linear connections
on the quantum group $GL_{q}(n)$ and study their properties.

\bigskip
\noindent
{\bf Proposition 5.1:\ }
{\it Let $\sigma$ be any generalized permutation.
The map $\nabla_{0}^\sigma$  defined by
\begin{eqnarray}
 &\nabla_{0}^\sigma: \Omega^{1}\rightarrow \Omega^{1}\otimes_{\cal A}\Omega^{1}
&
\nonumber\\
\label{deflincon}
& \nabla_{0}^\sigma(\omega)=\theta\otimes\omega-\sigma(\omega\otimes\theta). &
\end{eqnarray}
is a linear connection associated to $\sigma$.}

\noindent
{\it Proof. } Calculate first
\begin{equation}
\nabla_{0}^\sigma(a\omega)=([\theta,a]+a\theta)\otimes\omega-
\sigma(a\omega\otimes\theta)
\end{equation}
and then use the expression of the exterior derivative
and the bimodule property to obtain
\begin{equation}
\nabla_{0}^\sigma(a\omega)=da\otimes\omega+a\nabla_{0}^\sigma\omega.
\end{equation}
Similarly, calculate
\begin{eqnarray}
\nabla_{0}^\sigma(\omega a) &=& \theta\otimes \omega a+
\sigma(\omega\otimes([\theta,a]-\theta a)),\nonumber\\
 &=& \sigma(\omega\otimes da)+(\nabla_{0}^\sigma\omega)a.
\end{eqnarray}
This completes the Proof. $\clubsuit$
\bigskip

\noindent
{\bf Remarks:}

1-The linear connection $\nabla_{0}^\sigma$ can
be defined on any  differential calculus where the
exterior derivative is a graded commutator.
See \cite{mmm} for another example.

2-For $\sigma=-1$ the resulting covariant derivative $\nabla_{0}^\sigma$ is
$i\circ d$, where $i$ is the embedding of $\Omega^{2}$ into
$\Omega^{1}\otimes_{\cal A}\Omega^{1}$, by equation (2.21).

\bigskip

\noindent
{\bf Proposition 5.2:}
{\it The extension of $\nabla_{0}^{\sigma}$ to the tensor product
of $s$ copies of $\Omega^{1}$ is given by
\begin{equation}
\nabla_{0}^{\sigma} \ \nu\ =\ \theta\otimes\nu
+\sigma_{s}(\nu\otimes\theta),\
\forall\nu\in\Omega^{1}\otimes_{\cal A}\dots\Omega^{1}
\end{equation}
}

\noindent
{\it Proof.} A direct application of Proposition 3.3. $\clubsuit$

\bigskip

\noindent
{\bf Proposition 5.3: }{\it
There are no non-vanishing bimodule homomorphisms
from $\Omega^{1}$ to $\Omega^{1}\otimes_{\cal A}\Omega^{1}$.
}

\noindent
{\it Proof. } We will use the following Lemma proved in
\cite{fad,schu,zum}

\noindent
{\bf Lemma 5.4: } {\it
Let $c$ be the q-determinant of $T$,
\begin{equation}
c=det_{q}T=\sum_{p}(-q)^{l(p)}T^{1}_{p(1)}T^{2}_{p(2)}
\dots T^{n}_{p(n)},
\end{equation}
where the sum is over all permutations
on $n$ elements and $l(p)$ is the number
of transpositions in the permutation $p$.
 Then $c$ is in the center of $\cal A$ and verifies
$\omega c=q^{-2}c\omega$ for all $\omega$ in $\Omega^{1}$.}

An immediate consequence of the preceding Lemma
is the

\medskip
\noindent
{\bf Corollary 5.5: } $\nu c=q^{-4}c\nu,\ \forall  \nu
\in \Omega^{1}\otimes_{\cal A}\Omega^{1}$.

\medskip
We are now in position  to prove the Proposition.
Let $\phi$ be a bimodule homomorphism from
$\Omega^{1}$ to $\Omega^{1}\otimes_{\cal A}\Omega^{1}$ and
let $\omega \in \Omega^{1}$. By the homomorphism property
and the  Lemma 5.4, we get
\begin{equation}
\phi(\omega)c=q^{-2}c\phi(\omega).
\end{equation}
On the other hand, since $\phi(\omega)\in
\Omega^{1}\otimes_{\cal A}\Omega^{1}$,
by Corollary 5.5
 we obtain
\begin{equation}
\phi(\omega)c=q^{-4}c\phi(\omega).
\end{equation}
Comparing these two  equations we prove
the Proposition. $\clubsuit$

\bigskip
As a direct consequence of the preceding and of
the third remark following Definition 3.2 we obtain

\noindent
{\bf Theorem 5.6: }
{\it For any generalized permutation $\sigma$ on $GL_q(n)$,
there exists one and only one
associated linear connection, given by (\ref{deflincon}).
}

\bigskip
We now turn to the study
of some of the properties of the linear connection
$\nabla_{0}^\sigma$.

\bigskip
\noindent
{\bf Proposition 5.7: }
{\it For any generalized permutation $\sigma$,
the linear connection $\nabla_{0}^\sigma$ has vanishing torsion.}

\noindent
{\it Proof. } Calculate $\pi\circ\nabla_{0}^\sigma$
\begin{equation}
\pi\circ\nabla_{0}^\sigma \omega=\theta\wedge\omega+\omega\wedge\theta
=d\omega,
\end{equation}
where we have used the property (\ref{defsigma}).
The proof of the Proposition follows from (\ref{tor}).
$\clubsuit$

\bigskip
\noindent
{\bf Proposition 5.8: }
{\it For any generalized permutation $\sigma$,
the linear connection
$\nabla_{0}^\sigma$ has the expression
\begin{eqnarray}
\label{nablaleftinv}
\nabla_{0}^\sigma \omega^a &=&
(\Lambda - \sigma)\omega^a \otimes \theta \\
&=& \Big((\lambda_{1}-1)\Pi_{1}+(\lambda_{2}-1)\Pi_{2}+
(q^{2}-1)\Pi_{3}+(q^{-2}-1)\Pi_{4}\Big)\omega^{a}\otimes\theta,\nonumber
\end{eqnarray}
on the left invariant 1-forms $\omega^{a}$,
and
\begin{eqnarray}
\label{nablarightinv}
\nabla_0^\sigma \eta^a &=& (\Lambda^{-1} - \sigma)\eta^a \otimes \theta\\
&=& \Big((\lambda_{1}-1)\Pi_{1}+(\lambda_{2}-1)\Pi_{2}+
(q^{-2}-1)\Pi_{3}+(q^{2}-1)\Pi_{4}\Big)\eta^{a}\otimes\theta.\nonumber
\end{eqnarray}
on the right-invariant 1-forms $\eta^{a}$.}

\noindent
{\it Proof. } This is an immediate consequence of
the definition of $\Lambda$, the right
invariance of $\theta$ and Equation (\ref{fi}). $\clubsuit$

\bigskip
\noindent
{\bf Definition 5.9:} {\it A bicovariant linear connection, $\nabla$,
is defined  by the properties
\begin{eqnarray}
\label{leftcovariance}
 (1 \otimes\nabla)\Delta_{L}
&=&\Delta_{L} \nabla \ \ (\mbox{left covariance}), \\
\label{rightcovariance}
(\nabla\otimes 1)\Delta_{R}
 &=& \Delta_{R} \nabla\ \ (\mbox{right  covariance}).
\end{eqnarray}
}

\bigskip
\noindent
{\bf Proposition 5.10: }{\it
The linear connections associated to the generalized permutations
$\sigma_{\lambda_{1},\lambda_{2}}$ of formula
(\ref{fi}) are bicovariant.
}

\medskip
\noindent
{\it Proof. }
First, one sees that $\Lambda$ and $\sigma_{\lambda_1, \lambda_2}$ are
bicovariant. Then, using formula (\ref{nablaleftinv}) and the left
invariance of $\omega^a$ one sees that formula (\ref{leftcovariance}) is
true when applied to $\omega^a$.
Now, the 1-forms $\omega^a$ form  a basis of the  left
module  $\Omega^1$. Then, formula (\ref{defsigmaleft}) and the
previous result show that the associated
linear connection is left invariant.

For the right invariance, one has to consider the right  invariant
1-forms $\eta^a$, which constitute a basis of $\Omega^1$ as a
right module and
 formulas (\ref{nablarightinv}), (\ref{defsigmaright}). $\clubsuit$

\bigskip
The following Proposition allows one to calculate explicitly
the covariant derivative associated to a generalized permutation
given by equation (\ref{fi}):

\medskip

\noindent
{\bf Proposition 5.11: }
{\it Define $\nu,\gamma$ and $\beta$ by
\begin{equation}
\nu=q+q^{-1},\quad \gamma={{\lambda_{1}-\lambda_{2}}\over{q^{-2}-q^{2}}},\quad
\beta={{\lambda_{1}q^{2}-\lambda_{2}q^{-2}}\over{q^{-2}-q^{2}}};
\end{equation}
the linear connection associated to the generalized permutation,
$\sigma_{\lambda_{1},\lambda_{2}}$, acts on left-invariant
1-forms as follows:
\begin{eqnarray}
\displaystyle
\nabla^{\sigma_{\lambda_{1},\lambda_{2}}}_{0} \omega^{i}_{j}
=&-{{1}\over{\nu^2}}(1-\gamma-\beta)\ \omega^{i}_{k}\wedge\omega^{k}_{j}
-\gamma \ \omega^{i}_{k}\otimes\omega^{k}_{j} \nonumber\\
&+{1\over{2}}(1-\gamma+\beta)\ (\omega^{i}_{j}\otimes\theta+
\theta\otimes\omega^{i}_{j})\nonumber\\
&+{\lambda^{2} \over{2\nu^{2}}}
(1-\gamma-\beta)\
(\omega^{i}_{j}\otimes\theta-\theta\otimes\omega^{i}_{j}).
\end{eqnarray}
}

\noindent
{\it Proof.}
 First we note that $\sigma_{\lambda_{1},\lambda_{2}}$ of
(\ref{fi}), can be written as
\begin{equation}
\sigma_{\lambda_{1},\lambda_{2}}=(\lambda_{1}+1)\Pi_{1}+
(\lambda_{2}+1)\Pi_{2}-1,
\end{equation}
so that $ \nabla \omega^{i}_{j}$ can be expressed as
\begin{equation}
\nabla \omega^{i}_{j} =\theta\otimes\omega^{i}_{j}
+\omega^{i}_{j}\otimes\theta-\Big[(\lambda_{1}+1)\Pi_{1}+(\lambda_{2}
+1)\Pi_{2}\Big]\ \omega^{i}_{j}\otimes\theta.\label{int}
\end{equation}
It remains to calculate the term in the brackets of
(\ref{int}). We will do so by calculating
it for two different values of the couple
$(\lambda_{1},\lambda_{2})$ with the aid of
the  following two Lemmatae

\noindent
{\bf Lemma 5.12:} {\it The covariant
derivative associated to $\sigma_{\Lambda}$
acts on left-invariant 1-forms as follows:
\begin{equation}
\nabla^{\sigma_{\Lambda}}_{0}\omega^{i}_{j}
=-{{2}\over{\nu^{2}}}\omega^{i}_{k}\wedge\omega^{k}_{j}-
{\lambda^{2}\over{\nu^{2}}}
(\theta\otimes\omega^{i}_{j}-\omega^{i}_{j}\otimes\theta).
\end{equation}
}

\noindent
{\it Proof.} The Proof is a straigntforward calculation
exploiting the fact that $\sigma_{\Lambda}$ can be expressed in terms
of $\Lambda$ as
\begin{equation}
\sigma_{\Lambda}=-1+2{{(\Lambda+q^2)(\Lambda+q^{-2})} \over
{\nu^{2}}},
\end{equation}
as well as the  equation
\begin{equation}
d\omega^{i}_{j}=(1-\Lambda)\theta\otimes\omega^{i}_{j}=
-\omega^{i}_{k}\wedge\omega^{k}_{j},
\end{equation}
which allows to eliminate $\Lambda(\theta\otimes\omega^{i}_{j})$
in $\nabla^{\sigma_{\Lambda}}\omega^{i}_{j}$. $\clubsuit$

\medskip
\noindent
{\bf Lemma 5.13: }
{\it The covariant derivative associated to $\sigma_{R}$
is determined by
\begin{equation}
\nabla^{\sigma_{R}} \ dT^{i}_{j}=0.
\end{equation}
}

\noindent
{\it Proof.} Calculate the covariant derivative associated
to $\sigma_{R}$ of the two sides of equation (\ref{commutdt}). $\clubsuit$

The Proof of the Proposition is completed after  expressing
 $\nabla^{\sigma_{\lambda_{1},\lambda_{2}}} \omega^{i}_{j}$ in
terms of $\nabla^{{\sigma_{\Lambda}}}\omega^{i}_{j}$
and $\nabla^{{\sigma_{R}}}\omega^{i}_{j}$ as
\begin{eqnarray}
\nabla^{\sigma_{\lambda_{1},\lambda_{2}}}\omega^{i}_{j}
=&{1 \over{2}}(1-\gamma+\beta)(\theta\otimes\omega^{i}_{j}+
\omega^{i}_{j}\otimes\theta) \nonumber\\
&+{1 \over{2}}(1-\gamma-\beta)\nabla^{\sigma_{\Lambda}}\omega^{i}_{j}
+\gamma\nabla^{\sigma_{R}}\omega^{i}_{j}.
\end{eqnarray}
This equation is obtained after the evaluation of
$\Pi_{k}\ \omega^{i}_{j}\otimes\theta,\ k=1,2$ in terms of
$\nabla^{\sigma_{\Lambda}}\omega^{i}_{j}$ and
$\nabla^{\sigma_{R}}\omega^{i}_{j}$. $\clubsuit$

\bigskip

Finally, we consider the limit of
the linear connections determined above
when $q\rightarrow 1$.  In this limit the differential
calculus tends to the usual commutative differential
calculus. The 1-form $\theta$ has a singular
limit but $\lambda\theta$ tends to the right and left invariant
1-form $\alpha$ on $GL(n)$.
First of all, a necessary condition
for the limit to be non-singular is that the generalized permutation
tend to the flip operator that is
$\lambda_{1}\rightarrow 1$
and  $\lambda_{2}\rightarrow 1$.
A more precise statement, giving a necessary and
sufficient condition for the limit to be non-singular is the following:

\noindent
{\bf Proposition 5.14:} {\it Let
\begin{equation}
\mu_{i}={{\lambda_{i}-1} \over {\lambda}},\quad  i=1,2,
\end{equation}
the linear connection associated to $\sigma_{\lambda_{1},\lambda_{2}}$
admits a non-singular limit iff $\mu_{1}$ and $\mu_{2}$ have  finite
limits $\mu_{i}|_{q=1}$ when $q$ tends to 1.
The linear connection, in the limit, is
determined by
\begin{eqnarray}
\displaystyle
\nabla \omega^{i}_{j}=
&-{1 \over{2}} (1-\gamma_{0})\ \omega^{i}_{k}\wedge\omega^{k}_{j}-
\gamma_{0}\ \omega^{i}_{k}\otimes\omega^{k}_{j},\nonumber\\
\displaystyle
&-{\mu_{0} \over{2}}\ (\alpha\otimes\omega^{i}_{j}+
\omega^{i}_{j}\otimes\alpha),
\end{eqnarray}
where
\begin{equation}
\gamma_{0}={{\mu_{2}|_{q=1}-\mu_{1}|_{q=1}} \over {2}},\quad
\mu_{0}={{\mu_{2}|_{q=1}+\mu_{1}|_{q=1}} \over {2}}.
\end{equation}
}

\noindent
{\it Proof. } A direct application of Proposition
5.11. $\clubsuit$

\noindent
{\bf Remark:} When $\mu_{1}$ and $\mu_{2}$ tend to 0, which
is the case of $\sigma_{\Lambda}$,
$\gamma_{0}$
and $\mu_{0}$ vanish and the limiting linear connection is given by
\begin{equation}
\nabla \omega^{i}_{j}=
-{1 \over{2}}  \omega^{i}_{k}\wedge\omega^{k}_{j}.
\end{equation}

\section{Conclusion}

The main result of this paper is
the existence and uniqueness, for generic $q$,
of the linear connection associated to
a given generalized permutation.
This connection is  bicovariant
and torsion-free. This is in contrast to the
commutative case ($q=1$) where there are an
infinite number of linear connections not
necessarily bicovariant and torsion-free
and where the generalized
permutation is constrained to be the flip operator. It is
also in contrast to the cases with $q$ a root of unity
where Proposition (5.3) is not in general valid.
The arbitrariness in the deformed case lies merely in
the generalized permutation for which we
 have found a two parameter
family  (equation {\ref{fi}).
These parameters may be arbitrary functions of $q$
and are constrained by the involution property (Proposition 4.3).
The commutative limit is non-singular for a class of such functions which
tend to the identity when $q\rightarrow 1$.
The commutative limit of the linear connection
is a subset of right and left invariant
linear connections on $GL(n)$.

We have used the differential calculus
\ref{cal} to obtain our results. Had we used
another differential calculus with the usual commutative
limit the qualitative aspects of our conclusions, in particular
the uniqueness of the linear connection associated to a given
generalized permutation,
are expected to remain the same.


\begin{thebibliography}{99}
\bibitem{asc} {Aschieri P., Castellani L., Int. Jour. Mod. Phys. A8
(1993) 1667.}

\bibitem{ber} {Bernard D., Prog. Theoretical Phys. Supp. 102 (1990) 49.}

\bibitem{car} {Carow-Watamura U., Schlieker M., Watamura S., Weich W.,
Commun. Math. Phys. 142 (1991) 605.}

\bibitem{cham} {Chamseddine A.H., Felder G., Fr\"ohlich J., Commun.
Math. Phys. 155 (1993) 205.}

\bibitem{con} {Connes A., {\it Non-commutative differential
geometry}, Inst. Hautes Etudes Sci. Publ. Math. 62 (1986).}

\bibitem{dri} {Drinfeld V.G., Dokl. Akad. Nauk SSSR 283 (1985)
1060-1064; {\it Quantum Groups} Proc. Int. Congr. Math. (Berkley), vol 1,
Academic Press, New York, 1986 pp. 798-820.}

\bibitem{dub} {Dubois-Violette M., Madore J., Masson T., Mourad J.,
Lett. Math. Phys. (to appear)}

\bibitem{dub1} {Dubois-Violette M., Michor P.,{\it
Connections on central bimodules}, { Preprint LPTHE 94/100}}

\bibitem{fadd}{Faddeev L.D., Pyatov P.N., {\it The differential
calculus on quantum linear groups},  Preprint, hep-th 9402070.}

\bibitem{fad} {Faddeev L.D., Reshetikhin N.Yu., Takhtajan L.A.,
Algebra i Analiz 1 (1989) 178 (English translation:
Leningrad Math. J. 1
(1990) 193.)}

\bibitem{jim} {Jimbo M., Lett. Math. Phys. 11 (1986) 247-252;
Commun. Math. Phys. 102 (1986) 537.}

\bibitem{jur} {Jurco B., Lett. Math. Phys. 22 (1991) 177.}

\bibitem{kal} {Kalau W., Walze M., {\it Gravity, non-commutative
geometry and the Wodzicki residue}, Preprint MZ-TH/93-38.}

\bibitem{kas}{Kastler D., {\it The Dirac operator and gravitation}, Preprint
CPT 93/P.2970.}

\bibitem{keh} {Kehagias A., Madore J., Mourad J., Zoupanos G.,
{\it Linear connections on extended space-time}, J. Math. Phys. (to appear).}

\bibitem{kul} {Kulish P.P., Reshetikhin N. Yu., Zap. Nauchn.
Sem. Leningrad. Otdel. Mat. Inst. Steklov. (LOMI) 101
(1981) 101.}

\bibitem{mad2} {Madore J., Phys. Rev. D 41 (1990) 3709.}

\bibitem{mad1} {Madore J., {\it Introduction to non-commutative geometry and
its physical applications},
Cambridge University Press (1995)}

\bibitem{mmm} {Madore J., Masson T., Mourad J., Class. Quant. Grav.
 { 12} (1995) 1429. }

\bibitem{maj} {Majid S., Int. Jour. Mod. Phys. A5 (1990) 4689.}

\bibitem{mal} {Maltsiniotis G., C.R. Acad. Sci. Paris, 331 (1990)
831.}

\bibitem{man} {Manin Yu. I., {\it Quantum groups and
non-commutative geometry}, Preprint CRM-1561,
Montr\'eal, 1988.}

\bibitem{mou} {Mourad J., Class. Quant. Grav. 12 (1995) 965.}

\bibitem{sch} {Schmudgen K., Schuler A.,  Commun. Math. Phys.
167 (1995) 635.}

\bibitem{schu} {Schupp P., Watts P., Zumino B., Lett.
Math. Phys. 25 (1992) 139.}

\bibitem{skl} {Sklyanin E.K., Takhtadzhyan L.A., Faddeev L.D.,
I. Theoret. Mat. Fiz. 40 (1979) 194.}

\bibitem{schi} {Schirrmacher A., {\it Remarks on use of R-matrices,} in
Proc. of the first Max Born Symposium, (R Gielerak {\it et al} eds.),
Kluwer Academic Publishers (1992) p.55.}

\bibitem{sud} {Sudbery A., Phys. Lett. B284 (1992) 61.}

\bibitem{wes} {Wess J., Zumino B., Nucl. Phys. (Proc. Suppl.) B18 (1990) 302.}

\bibitem{woro} {Woronowicz S.L., Commun. Math. Phys. 111 (1987) 613.}

\bibitem{wor} {Woronowicz S.L., Commun. Math. Phys. 122 (1989) 125.}

\bibitem{zum} {Zumino B., {\it Differential calculus on quantum spaces
and quantum groups,} Preprint LBL-33249, UCB-PTH-92/41.}

\end{thebibliography}
\end{document}